\documentclass[reqno, 11pt]{amsart}
\usepackage{color}
\usepackage{graphicx} 
\usepackage{hyperref}

\newtheorem{theorem}{Theorem}[section]
\newtheorem{thm}[theorem]{Theorem}
\newtheorem{cor}[theorem]{Corollary}
\newtheorem{lem}[theorem]{Lemma}

\newtheorem{rem}[theorem]{Remark}

\theoremstyle{definition} \theoremstyle{remark}
\numberwithin{equation}{section}

\newcommand{\C}{\mathbb C}
\newcommand{\R}{\mathbb R}

\renewcommand{\P}{\mathbb P}
\newcommand{\N}{\mathbb N}
\newcommand{\g}{\mathfrak g}
\renewcommand{\k}{\mathfrak k}
\renewcommand{\a}{\mathfrak a}

\newcommand{\n}{\mathfrak n}
\newcommand{\p}{\mathfrak p}

\newcommand{\GL}{\rm GL}
\newcommand{\mtx}[1]{\begin{bmatrix} #1 \end{bmatrix}} 

\newcommand{\diag}{\operatorname{diag}}
\newcommand{\ad}{\operatorname{ad}}
\newcommand{\Ad}{\operatorname{Ad}}
\newcommand{\End}{\operatorname{End}}
\newcommand{\Aut}{\operatorname{Aut}}

\newcommand{\U}{\operatorname {U}}
\newcommand{\im}{\operatorname{Im}}

\newcommand{\be}[2]{\begin{#1} #2 \end{#1}}  

\usepackage{xcolor, soul}

\title{Extensions of Yamamoto-Nayak's Theorem}
\author{Huajun Huang}
\address{Auburn University, Alabama, USA}
\email{huanghu@auburn.edu}
\author{Tin-Yau Tam}
\address{University of Nevada, Reno, USA}
\email{ttam@unr.edu}
\date{\today}

\begin{document}
\begin{abstract}
 A result of Nayak asserts that $\underset{m\to
\infty}\lim |A^m|^{1/m}$ exists for each $n\times n$ complex matrix $A$, where $|A| = (A^*A)^{1/2}$, and the limit is given in terms of the spectral decomposition. We extend the result of Nayak, namely, we prove that the limit 
of $\underset{m\to
\infty}\lim |BA^mC|^{1/m}$ exists for any $n\times n$ complex matrices $A$, $B$, and  $C$, where $B$ and $C$ are nonsingular; the limit is obtained and is independent of $B$. We then provide generalization in the context of real semisimple Lie
groups. 
\end{abstract}


\maketitle

\section{Introduction}

Let $\N$  (resp. $\R$, $\C$) be the set of positive integers (resp.  real numbers, complex numbers).
Let $M_n(\C)$ (resp. $M_n(\R)$) denote the set of $n\times n$ complex (resp. real) matrices,  $\GL_n(\C)$ the  group of $n\times n$ complex nonsingular matrices, and $\U(n)$ the  group of $n\times n$ unitary matrices. 
Let $\P_n$ (resp. $\overline{\P}_n$) be the set of positive definite (resp. positive semidefinite) matrices in $M_n(\C)$.  
For $A\in M_n(\C)$, let $|A|=(A^*A)^{1/2}$ and let $\lambda_1(A),\ldots,\lambda_n(A)$ denote the eigenvalues of $A$, counting multiplicities, in the way that
$|\lambda_1(A)|\ge |\lambda_2(A)|\ge \cdots\ge|\lambda_n(A)|$. Define the following:
\begin{itemize}
\item $\lambda(A)=(\lambda_1(A),\ldots,\lambda_n(A))$ the $n$-tuple of eigenvalues;
\item $|\lambda|(A)=(|\lambda_1 (A)|,\ldots,|\lambda_n (A)|)$ the $n$-tuple of eigenvalue moduli with non-increasing order;
\item $s(A)=(s_1(A),\ldots,s_n(A))$ the $n$-tuple of singular values of $A$ with $s_1(A)\ge s_2(A)\ge \cdots\ge s_n(A)$.
\end{itemize}

 Beruling-Gelfand's spectral radius formula asserts that \cite [p.70]{S} for $A\in M_n(\C)$: 
\begin{equation}\label{spectral}
 \lim_{m\to \infty} \|A^m\|^{1/m} = \rho(A),
 \end{equation} where
$\rho(A)$ denotes the spectral radius of $A$ and $\|A\|=s_1(A)$ denotes the
spectral norm of $A$.   Indeed it is true for all matrix norms \cite [p.349]{HJ13}.  It is an interesting asymptotic result  as it relates algebraic and analytic properties of $A$ in a nice way. 
  Yamamoto generalized the result to all singular values  \cite{Y} 
\begin{equation}\label{yamamoto}
 \lim_{m\to \infty}
[s_i(A^m)]^{1/m} = |\lambda_i(A)|,\quad i=1, \dots , n.
\end{equation}
 Tam and Huang extended
\eqref{yamamoto}
 in the context of semisimple Lie groups \cite{TamH06}.
 They also proved  that \cite[Theorem 1.2]{HuangTam08} for $A\in M_n(\C)$ and nonsingular $B, C\in M_n(\C)$, 
\begin{equation}
    \label{BA^mC SV limit}
\lim_{m\to\infty} s(BA^mC)^{1/m}=|\lambda|(A).
\end{equation} 
Denote by $|A| = (A^*A)^{1/2}$, the positive semidefinite component of the polar decomposition of 
\begin{equation}\label{polar}
A = U|A|,\qquad U\in \U(n)
\end{equation}
 So $$s_i(A) =s_i(|A)|= \lambda_i(|A|),\quad i=1, \dots, n.$$
Yamamoto's result is about scalar convergence but not the convergence of the matrix sequence $\{|A^m|^{1/m}\}_{m\in \N}$. From Yamamoto's result, the limit should be in the unitary orbit of $|A|$ via similarity, if the  matrix sequence converges. Very recently,  Nayak  \cite{Nayak} has proved that the matrix sequence does converge.

\begin{theorem}[Nayak, 2023] \label{Nayak}  Let $A \in M_n(\C)$ and $\{\gamma_1,  \dots  , \gamma_s\}$ be the set of eigenvalue moduli  of $A$ such that $\gamma_1>\gamma_2>\cdots>\gamma_s\ge 0$.
 Let $A =  S +N$ be the  
 (additive)  Jordan-Chevalley decomposition of $A$ into its commuting diagonalizable  $S$  and nilpotent $N$. 
For $1\leq j\leq s$, let $E_j$ be the orthogonal projection onto the subspace of $\C^n$ 
spanned by the eigenvectors of   $S$  corresponding to eigenvalues with moduli no greater than $\gamma_j$ and set $E_{s+1}=0$. Then the following assertions hold:
\begin{itemize}
\item[(i)] The sequence $\{|A^m|^{1/m}\}_{m\in \N}$ converges to the   positive semidefinite  matrix $\sum_{i=1}^s \gamma_i(E_i-E_{i+1})$. 
\item[(ii)] 
Let $\im E_j$ denote  the image (that is, the column space)  of $E_j$.
A nonzero vector $x\in \C^n$ is in $\im E_j\setminus \im E_{j+1}$ if and only if $$\lim _{m\to \infty}\|A^m x\|^{1/m} =\gamma _j.$$  
\item [(iii)] The set $\im E_j\setminus \im E_{j+1}$ is invariant under the action of $A^k$ for every $k \in \N$.
\end{itemize}
\end{theorem}

The theorem is stated using the language of linear transformation. In this paper we will establish  an extension of Nayak's result in matrix  version and  our proof is  different from Nayak's. We further  give an extension of Theorem \ref{Nayak}(i)   in the context of real semisimple Lie group. 


\section{The limit  of $|BA^mC|^{1/m}$}

For $n\times n$ Hermitian matrices $A$ and $B$, we use $A\ge B$ to denote the  L\"owner  order, that is,
$A\ge B$ if and only if $A-B$ is positive semidefinite.

The following result is analogous to the Sandwich Theorem. 

\begin{lem}\label{thm: Sandwich}
If $\{A_m \}_{m\in\N}$, $\{ B_m \}_{m\in\N}$, and $\{C_m \}_{m\in\N}$ are sequences  in $\P_n$ such that
\begin{enumerate}
\item $A_m\le B_m\le C_m$ for $m\in\N$, and
\item $\lim_{m\to\infty} A_m=B= \lim_{m\to\infty} C_m$ for some $B\in\P_n$,
\end{enumerate} 
then $\lim_{m\to\infty} B_m=B$.
\end{lem}

\begin{proof}
Since $\lim_{m\to\infty} C_m=B$, the sequence  $\{C_m \}_{m\in\N}$ is contained in a compact set of $\P_n$. 
So is the sequence  $\{B_m \}_{m\in\N}$. For every limit point $B'$ of $\{B_m \}_{m\in\N}$ we have
$B\le B'\le B$ so that $B'=B$. Therefore, $\lim_{m\to\infty} B_m=B$.
\end{proof}

\begin{lem}\label{Thm: B_mA_m limit}
Suppose that the sequence $\{ B_m \}_{m\in\N} \subseteq\GL_n(\C)$  satisfies
\be{equation}{\label{SV limit bound}
\lim_{m\to\infty} s_n(B_m)^{1/m}=1=\lim_{m\to\infty} s_1(B_m)^{1/m}.
}
Then
for any $\{A_m \}_{m\in\N}\subseteq M_n(\C)$ and $A\in\P_n$,
\be{equation}{
\lim_{m\to\infty} |A_m|^{1/m} =A\quad \Rightarrow\quad
\lim_{m\to\infty} |B_mA_m|^{1/m}=A.
}
\end{lem}

\begin{proof}
In the  L\"owner  order,
\be{equation}{
s_n(B_m)^2  I_n\le B_m^*B_m \le  s_1(B_m)^2 I_n.
}
Hence for $A_m\in M_n(\C)$:
\be{eqnarray}{
A_m^*(B_m^*B_m-s_n(B_m)^2  I_n)A_m &\ge& 0,
\\
A_m^*(s_1(B_m)^2  I_n-B_m^*B_m)A_m &\ge& 0,
}
which give 
\be{equation}{
s_n(B_m)^2 A_m^* A_m\le A_m^*  B_m^*B_m A_m\le s_1(B_m)^2   A_m^* A_m.
}
By \cite[Theorem 1.5.9]{Bhatia07}, the function $X\mapsto X^{1/(2m)}$ is monotone on $\P_n$. So 
\be{equation}{
s_n(B_m)^{1/m} |A_m|^{1/m} \le |B_m A_m|^{1/m} \le s_1(B_m)^{1/m} |A_m|^{1/m}.
} 
By \eqref{SV limit bound} and Lemma \ref{thm: Sandwich}, when
$\lim_{m\to\infty}|A_m|^{1/m}=A$, we get 
$$\qquad \qquad \qquad\lim_{m\to\infty}|B_mA_m|^{1/m}=\lim_{m\to\infty}|A_m|^{1/m}=A.
\qquad \qquad \qquad \qedhere$$
\end{proof}

\begin{rem}\rm
Condition \eqref{SV limit bound} is equivalent to each of the following conditions:
 \begin{enumerate}
     \item $\lim_{m\to\infty} |B_m|^{1/m}=I_n$.
     \item $\lim_{m\to\infty} \rho(B_m)^{1/m}=1=\lim_{m\to\infty} \rho(B_m^{-1})^{1/m}$. 
     \item For every norm $\|\cdot\|'$ of $M_n(\C)$, $$\lim_{m\to\infty} {\|B_m\|'}^{1/m}=1=\lim_{m\to\infty} {\|B_m^{-1}\|'}^{1/m},$$
     since every norm is equivalent to the spectral norm and $s_n(B_m)=\|B_m^{-1}\|^{-1}$.
 \end{enumerate}
\end{rem}

\begin{rem}\label{rectangle B_mA_m limit}\rm
Lemma \ref{Thm: B_mA_m limit} still holds when  $\{A_m \}_{m\in\N}$ are $n\times r$ complex matrices and $A\in\P_r$,
and the proof is analogous.
\end{rem}

Given $A\in M_n(\C)$ and $I, J\subseteq [n]$, let $A[I, J]$ denote the submatrix of  $A$ with rows indexed by $I$ and columns indexed by $J$, and
 abbreviate $A[p]=A[[p],[p]]$ for $p\in [n]$. 

\begin{lem}\label{thm: limit diag power low tri}
Let $D=\diag(d_1,\ldots, d_n)\in M_n(\R)$, where $d_1\ge\cdots\ge d_n\ge 0$. Then for any  nonsingular lower triangular matrix $L\in\GL_n(\C)$, we have
\be{equation}{\label{lim diag power}
\lim_{m\to\infty} | D^mL|^{1/m}=D.
}
\end{lem}

\begin{proof}
Rewrite 
\be{equation}{
D=\diag(d_1,\ldots, d_n)= \mu_1 I_{n_1}\oplus \cdots \oplus\mu_k I_{n_k},
}
where $k\in\N$, $\mu_1> \cdots>\mu_k\ge 0$, and $n_1,\ldots,n_k\in\N$ such that $n_1+\cdots+n_k=n$. If $k=1$ then $D=d_1I_n$ and \eqref{lim diag power} is obviously true. 
We assume $k\ge 2$ in the following proof.

Denote
\be{equation}{
| D^m L|^{1/m} := X_m:=U_mD_mU_m^*,
} 
where $U_m\in \U(n)$ and 
\be{equation}{
D_m  =\diag(d_{m,1},d_{m,2},\ldots,d_{m,n}),
}
with
$d_{m,1}\ge d_{m,2}\ge \cdots\ge d_{m,n}\ge 0$. 
The spectral norm $\|\cdot\|$ is submultiplicative, that is  $\|AB\|\le\|A\|\|B\|$ for all $A, B\in M_n(\R)$. So
\be{eqnarray}{
\| X_m\| =\|D^m L\|^{1/m}
\le ( \|D\|^m \|L\|)^{1/m}
= (\|L\|)^{1/m} d_1.
}
Since $\lim_{m\to\infty} (\|L\|)^{1/m}=1$, the sequence $\{X_m\}_{m\in\N}$ is contained in a bounded closed subset of $\overline{\P}_n$, which is compact. It remains to show that every limit point of $\{X_m\}_{m\in\N}$ equals to $D$.

For every $p\in [n]$ and $K\subseteq[n]$ with $|K|=p$, first
apply Cauchy-Binet formula to   $X_m^{2m}=U_mD_m^{2m}U_m^*$ to have
\be{eqnarray}{\notag
&&
\det((X_m^{2m})[[p],K]) 
\\ \notag
&=&  
\sum_{\substack{J\subseteq[n]\\|J|=p}}
 \det(U_m[[p],J])\; \det(D_m^{2m}[J,J])\; \det(U_m^*[J,K]) 
\\  \label{det power submatrix 1}
&=&
\sum_{\substack{J\subseteq[n]\\|J|=p}}
 \det(U_m[[p],J])\;   \overline{\det(U_m[K,J])}\; (\prod_{j\in J} d_{m,j})^{2m}.
}
 Applying  Cauchy-Binet formula to   $X_m^{2m}=|D^mL|^2=L^* D^{2m}L$  yields 
\be{eqnarray}{\notag
\det((X_m^{2m})[[p],K]) 
&=& \det( (L^* D^{2m} L) [[p],K])
\\\notag
&=&
\sum_{\substack{J\subseteq[n]\\|J|=p}} 
\det(L^*[[p],J])\det(D^{2m}[J,J])\det (L[J,K])
\\ \label{det power submatrix 2}
&=&
\sum_{\substack{J\subseteq[n]\\|J|=p}} 
 \overline{\det(L [J,[p]])}\det (L[J,K]) (\prod_{j\in J} d_j)^{2m}.  
}
In particular, \eqref{det power submatrix 1} and \eqref{det power submatrix 2} imply that for  $p\in [n]$,
\be{eqnarray}{
\label{power expression 1}
\det((X_m^{2m})[p]) &=&
\sum_{\substack{J\subseteq[n]\\|J|=p}}  |\det(U_m[[p],J])|^2 (\prod_{j\in J} d_{m,j})^{2m}
\\ \label{power expression 2}
&=& \sum_{\substack{J\subseteq[n]\\|J|=p}} |\det(L [J,[p]])|^2  (\prod_{j\in J} d_j)^{2m}
\\\label{power expression 3}
&\ge&
|\det(L  [p] )|^2  (\prod_{j=1}^{p} d_j)^{2m}.
}
We also have
\be{equation}{\label{det unitary principal minor}
\sum_{\substack{J\subseteq[n]\\|J|=p}}  |\det(U_m[[p],J])|^2=\det((U_mU_m^*)[p])=1.
}

Suppose that $X$ is a limit point of $\{X_m\}_{m\in\N}$, and  $\{X_{m_{\ell}}\}_{\ell\in\N}\subseteq \{X_m\}_{m\in\N}$   satisfying  
\be{equation}{
X =\lim_{\ell\to\infty} X_{m_{\ell}} =\lim_{\ell\to\infty}  U_{m_{\ell}}  D_{m_{\ell}} U_{m_{\ell}}^* .
}
By \eqref{BA^mC SV limit}, $s(X)=|\lambda|(D)=(d_1,\ldots,d_n).$ So $\lim_{\ell\to\infty} D_{m_{\ell}}=D$.
The sequence $\{U_{m_{\ell}}\}_{\ell\in\N}$ is in the compact unitary group  $\U(n)$; hence it has a converging subsequence.
By refining and reindexing the subsequence, we may assume that $\lim_{\ell\to\infty} U_{m_{\ell}}=U$ for  some  $U \in \U(n)$.
Then $X=UDU^*$. 

Since  $\lim_{m\to\infty} D_{m}=D=\mu_1I_{n_1}\oplus \cdots \oplus\mu_k I_{n_k}$, 
for $J\subseteq[n]$, $|J|=n_1$ and $J\ne [n_1]$:
\be{equation}{
\lim_{m\to\infty}   \Big(\prod_{j\in J}  d_{m,j}\Big )^{2m}/\Big(\prod_{i\in [n_1]}d_i\Big)^{2m}=
\lim_{m\to\infty}  \Big (\prod_{j\in J} \frac{d_{m,j}}{\mu_1} \Big)^{2m}=0.
}
Consider $\{X_{m_{\ell}}\}_{\ell\in\N}$.
By \eqref{power expression 1}, \eqref{power expression 2}, and \eqref{det unitary principal minor}  for $p=n_1$, we have
\be{eqnarray}{\notag
&&\lim_{\ell\to\infty} \det((X_{m_{\ell}}^{2m_{\ell}})[n_1])/\Big(\prod_{i\in [n_1]}d_i\Big)^{2m_{\ell}}\\
&=& |\det(U  [n_1] )|^2 \lim_{\ell\to\infty} \Big (\prod_{j\in [n_1]} \frac{d_{m_{\ell},j}}{\mu_1}\Big)^{2m_{\ell}} 
\\ \label{limit diag block}
&=& |\det(L  [p] )|^2 >0.
}
 Now let $K\subseteq [n]$, $|K|=n_1$, and $K\ne [n_1]$. Since $L$ is lower triangular, we have
$\det (L[[n_1],K])=0$, so that by   \eqref{det power submatrix 2},
\be{equation}{
\lim_{\ell\to\infty} \det((X_{m_{\ell}}^{2m_{\ell}})[[n_1],K])/\Big(\prod_{i\in [n_1]}d_i\Big)^{2m_{\ell}}=0.
}
Together with \eqref{det power submatrix 1}, we have
\be{eqnarray}{\notag
0&=&\lim_{\ell\to\infty} \det((X_{m_{\ell}}^{2m_{\ell}})[[n_1],K])/\Big(\prod_{i\in [n_1]}d_i\Big)^{2m_{\ell}}
\\ \label{limit off diag block}
&=& \det(U [n_1])\; \overline{\det(U[K,[n_1]])}\;  \lim_{\ell\to\infty} \Big(\prod_{j\in [n_1]} \frac{d_{m_{\ell},j}}{\mu_1}\Big)^{2m_{\ell}}.
}
Therefore, \eqref{limit diag block} and \eqref{limit off diag block} imply that
$\det(U  [n_1] )\ne 0$ and $\det(U[K,[n_1]])=0$ for any $K\subseteq [n]$, $|K|=n_1$, and $K\ne [n_1]$. 
The rows of $U  [n_1]$ form a basis of the row vector space of dimension $n_1$.
If there exists $j\in [n]\setminus[n_1]$ such that $U[\{j\},[n_1]]\ne 0$, then 
there is  $(c_1,\ldots, c_{n_1})\ne (0,\ldots,0)$ such that
\be{equation}{
U[\{j\},[n_1]]=\sum_{i\in [n_1]} c_i U[\{i\},[n_1]].
}
Suppose $c_s\ne 0$ for certain $s\in [n_1]$. Then for $K=([n_1]\cup\{j\})\setminus\{s\}$ we get 
$\det(U[K,[n_1]])\ne 0$, which is a contradiction. Hence   $U\in \U(n)$  has the form
\be{equation}{
U=\mtx{U[n_1] &*\\0 &*}=\mtx{U[n_1] &0\\0 &*}.
}

For each $t\in [k-1]$, applying  
\eqref{det power submatrix 1} -- \eqref{det unitary principal minor} for $p=n_1+\cdots+n_{t}$  
and analogous arguments to $\{X_{m_{\ell}}\}_{\ell\in\N}$, we  have 
\be{equation}{
U= \mtx{U[n_1+\cdots+n_{t}] &0\\0 &*}.
}
Therefore, $U=U_1\oplus\cdots\oplus U_k$, where each $U_i\in \U(n_i)$. 

Overall,  
\be{equation}{
X=UDU^*=(U_1\oplus\cdots\oplus U_k)(\mu_1I_{n_1}\oplus \cdots \oplus\mu_k I_{n_k})(U_1^*\oplus\cdots\oplus U_k^*)=
D.
}
So every limit point of $\{X_m\}_{m\in\N}$ equals to $D$ and thus 
\eqref{lim diag power} is proved.
\end{proof}

The complete multiplicative Jordan decomposition (CMJD) of $A\in \GL_n(\C)$  is $A=EHU$, where
\begin{itemize}
\item
$E$ is {\em elliptic}, that is, $E$ is diagonalizable and $|\lambda|(E)=(1,\ldots,1)$;
\item
$H$ is {\em hyperbolic}, that is, $H$  is diagonalizable and $\lambda(H)=|\lambda|(H)$;
\item
$U$ is {\em unipotent}, that is, $\lambda(U)=(1,\ldots,1)$;
\item
the components $E$, $H$, and $U$  mutually  commute. 
\end{itemize}
The CMJD of $A$ is unique \cite {Helgason}. Explicitly, suppose that $A$ is similar to the Jordan canonical form 
\begin{equation}\label{JCF}
    A=M\big(\bigoplus_{i=1}^{k} J_{A}(\mu_i) \big)M^{-1},
\end{equation}  where
$\mu_1,\ldots,\mu_k$ are the distinct eigenvalues of $A$ and each $J_{A}(\mu_i)$ for $i\in[k]$ is the direct sum of Jordan blocks of $A$ associated to the eigenvalue $\mu_i$; we arrange $\mu_1,\ldots,\mu_k$ in the way that $|\mu_1|\ge \cdots\ge |\mu_k|$;
let $n_i$ $(i\in [k])$ be the algebraic multiplicity of the eigenvalue $\mu_i$. The CMJD $A=EHU$ is given by:
\begin{eqnarray}\label{elliptic}
E &=& M\big(\bigoplus_{i=1}^{k} \frac{\mu_i}{|\mu_i|}I_{n_i}\big) M^{-1},
\\ \label{hyperbolic}
H&=&M\big(\bigoplus_{i=1}^{k}  |\mu_i| I_{n_i} \big)M^{-1},
\\ \label{unipotent}
U &=& M\big(\bigoplus_{i=1}^{k} \frac{1}{ \mu_i }J_{A}(\mu_i)\big)M^{-1},
\end{eqnarray}
in which $\bigoplus_{i=1}^{k}  |\mu_i| I_{n_i}=\diag(|\lambda_1 (A)|,\ldots,|\lambda_n (A)|).$

When $A\in M_n(\C)$ is singular, the CMJD of $A$ is not well-defined. However, $A$ still has the hyperbolic component $H$ defined by \eqref{hyperbolic}, where we   assume that
$|\mu_1|\ge \cdots\ge |\mu_{k-1}|>|\mu_k|=0$. Let
\begin{eqnarray}
E' &=& M\big( \bigoplus_{i=1}^{k-1} \frac{\mu_i}{|\mu_i|}I_{n_i}\oplus I_{n_k}\big)M^{-1},
\\
U' &=& M\big[ \bigoplus_{i=1}^{k-1} \frac{1}{ \mu_i }J_{A}(\mu_i)\oplus (I_{n_k}+J_{A}(0))\big]M^{-1}.
\end{eqnarray}
 Then  $E'$ is elliptic, $U'$ is unipotent, $E'$, $H$, and $U'$  mutually  commute, and when $m$ is no less than the maximal size of the Jordan blocks of $A$ associated to eigenvalue $0$, we have
 \begin{equation}
 A^m={E'}^m H^m {U'}^m.
 \end{equation}

\begin{thm}\label{thm: Yamamoto extension}
Suppose $A\in M_n(\C)$. Let $H=MDM^{-1}$ be the hyperbolic element of $A$ given  in \eqref{hyperbolic}, in which 
$D=\diag(|\lambda_1 (A)|,\ldots,|\lambda_n (A)|).$
For any $B, C\in\GL_n(\C)$,  let 
$M^{-1}C=LQ$  in which $L$ is  lower triangular and $Q$ is   unitary. Then
\begin{equation}\label{BA^mC limit}
\lim_{m\to\infty} |BA^mC|^{1/m} = Q^* DQ,
\end{equation}
 in which the limit is independent of $B$  and the choice of $M$. Analogous result is true for $A\in M_n(\R)$, $B, C\in \GL_n(\R)$ and the unitary matrix $Q$ may be replaced by a real orthogonal matrix.
\end{thm}

\begin{proof}
We prove for the nonsingular case $A\in\GL_n(\C)$ and use the CMJD $A=EHU$. The proof for singular $A$ is similar, in which we use   $A^m={E'}^m H^m {U'}^m$ when $m$ is sufficiently large. 

For the nonsingular $A=EHU$, the components $E$, $H=MDM^{-1}$, and $U$ mutually  commute, so that 
\begin{eqnarray}\label{BA^mC decomp}
BA^mC 
= BE^m U^m MD^mM^{-1}C    
= (BE^m U^mM) (D^mLQ).
\end{eqnarray}
By Lemma \ref{thm: limit diag power low tri}, 
\begin{eqnarray} \notag
 \lim_{m\to\infty} |D^mLQ|^{1/m} 
 &=&
\lim_{m\to\infty} (Q^*L^*D^{2m}LQ)^{1/(2m)}
\\ \label{D^mLQ limit}
&=&
Q^*(\lim_{m\to\infty} |D^{ m}L|^{1/m})Q    
=Q^*DQ.
\end{eqnarray}
We claim that
\begin{equation}\label{BE^mU^mM}
 \lim_{m\to\infty} s_n(BE^m U^mM)^{1/m}=1 = \lim_{m\to\infty} s_1(BE^m U^mM)^{1/m}, 
\end{equation}
so that by \eqref{BA^mC decomp}, \eqref{D^mLQ limit}, and Lemma \ref{Thm: B_mA_m limit}, we get \eqref{BA^mC limit}. 
The independence of the choice of $B$ regarding  the limit \eqref{BA^mC limit}  is shown in  Lemma \ref{Thm: B_mA_m limit}. 
 Next we will   prove \eqref{BE^mU^mM}. 

By \eqref{elliptic}, the matrices $E^m$ and $E^{-m}$ for all $m\in\N$ are in the following compact subset of $\GL_n(\C)$: 
\begin{equation}\label{E compact set}
    \{M\diag(t_1,t_2,\cdots, t_n)M^{-1}\mid t_i\in \C,\ |t_i|=1,\ i\in [n]\}.
\end{equation} 
The function $A\mapsto \|A\|$ is continuous in $M_n(\C)$. So it has  a maximum  $c >0$ in the compact set \eqref{E compact set}. 
We get $\|E^m\|<c$ and $\|E^{-m}\|<c$ for all $m\in\N$. 

By \eqref{unipotent}, the unipotent element
$U=M(I+N)M^{-1}$ in which $N\in M_n(\C)$ is strictly upper triangular. Since $N^n=0$, for all $m\in\N$,
we have $U^m =M(I+N)^mM^{-1}$ and $U^{-m} =M(I+N)^{-m} M^{-1}$ in which the entries of matrices 
\begin{equation}
   (I+N)^m=\sum_{i=0}^{n-1} \binom{m}{i} N^i  
   \quad\text{and}\quad
   (I+N)^{-m}=\sum_{i=0}^{n-1} \binom{-m}{i} N^i 
\end{equation}
can be expressed as polynomials of $N$  of degrees less than $n$. So are the entries of $U^m$ and $U^{-m}$.
There exists a fixed  polynomial $g(x)\in\R[x]$ such that 
$\|U^m\|\le |g(m)|$ and 
$\|U^{-m}\|\le |g(m)|$ 
for all $m\in\N$. Therefore,
\begin{eqnarray}
    s_1(BE^m U^mM) &\le& \|B\|\|E^m\|\|U^m\|\|M\|\le (c  \|B\|\|M\|) |g(m)|,\qquad 
\\\notag 
 s_n(BE^m U^mM)
 &=& \|(BE^m U^mM)^{-1}\|^{-1} = \|M^{-1} U^{-m} E^{-m} B^{-1}\|^{-1} 
 \\ \notag &\ge & \|M^{-1}\|^{-1} \|U^{-m}\|^{-1} \|E^{-m}\|^{-1} \|B^{-1}\|^{-1} 
 \\
 &\ge &  (c ^{-1}s_n(B)s_n(M)) |g(m)|^{-1}.
\end{eqnarray}
Since $\lim_{m\to\infty} |g(m)|^{1/m}=1$ for every nonzero polynomial $g(x)$, 
the classical Sandwich Theorem implies \eqref{BE^mU^mM}. 

 The limit in \eqref{BA^mC limit} is described in terms of $Q$, which  depends on $M$.
We shall show that the limit in \eqref{BA^mC limit} is independent of the choice of $M$ 
as long as the hyperbolic element $H$ of $A$ satisfies that $H=MDM^{-1}$
for $D=\diag(|\lambda_1 (A)|,\ldots,|\lambda_n (A)|).$ 
 Suppose that $\hat M\in \GL_n(\C)$ is another choice, that is, $H=\hat MD{\hat M}^{-1}$. If $D = |\gamma_1| I_{m_1} \oplus \cdots \oplus |\gamma_s|I_{m_s}$, where
 $\gamma_1 > \cdots > \gamma_s \ge 0$, $m_1,\ldots,m_s\in \N$ and $m_1+\cdots+m_s=n$, 
 then $\hat M = M(M_1\oplus \cdots \oplus M_s)$, where $M_i\in \GL_{m_i}(\C)$. 
 Consider
 $$
 {\hat M}^{-1}C =(M_1^{-1}\oplus \cdots \oplus M_s^{-1})M^{-1}C =(M_1^{-1}\oplus \cdots \oplus M_s^{-1})LQ.
 $$
 Note that  $\hat L:=(M_1^{-1}\oplus \cdots \oplus M_s^{-1})L$ is in block lower triangular form and the main diagonal blocks are of size $m_1 , \dots, m_s$.
Performing Gram-Schmidt process on the rows of $\hat L$ from the top row to the bottom row, we have $\hat L = L_1\hat Q$ where $L_1$ is lower triangular and $\hat Q =  Q_1\oplus \cdots \oplus Q_s$ is unitary. 
 Hence  ${\hat M}^{-1}C =  L_1\hat Q Q$. Thus
 $$
({ \hat Q}Q)^{*} D({ \hat Q }Q)  = Q^{*}(\hat Q)^{*} D\hat Q Q =  Q^{*} DQ,
 $$
which is independent of the choice of $M$. 
\end{proof}

\begin{rem} \rm When $B=C = I_n$,  Theorem \ref{thm: Yamamoto extension} recovers   Nayak's result (see Theorem \ref{Nayak}). Suppose 
$$D=\diag(|\lambda_1 (A)|,\ldots,|\lambda_n (A)|)=\gamma_1 I_{m_1} \oplus \cdots \oplus \gamma_s I_{m_s},$$ where
 $\gamma_1 > \cdots > \gamma_s \ge 0$. 
 The hyperbolic part of $A$ is $H=MDM^{-1}$,
and
 $M=Q^*L^{-1}$ where $Q^*$ is unitary and $L^{-1}$ is lower triangular. 
 If   $M$ and $Q^*$ are partitioned   according to the column  partition $(m_1, \ldots,m_s)\vdash n$ 
 such that 
 $$M=[M_1\mid \cdots\mid M_s],\qquad Q^*=[Q_1\mid \cdots \mid Q_s],$$
then 
  $E_j$ in Theorem \ref{Nayak} is the orthogonal projection onto 
  $\im[M_j\mid\cdots\mid M_s]=\im [Q_j\mid \cdots\mid Q_s]$. The  columns of $Q^*$ are orthonormal. 
 Hence $E_j-E_{j+1}=Q_{j} Q_{j}^*$ and we get
   Theorem \ref{Nayak}(i): 
\begin{eqnarray}
   \lim_{m\to\infty} |A^m|^{1/m}
  =Q^*DQ= \sum_{j=1}^s \gamma_{j} Q _{j} Q_{j}^*=\sum_{j=1}^s \gamma_j(E_j-E_{j+1}).
\end{eqnarray}
From \eqref{JCF},  each $\im E_j=\im[M_j\mid\cdots\mid M_s]$ is the direct sum of generalized eigenspaces of $A$ with the eigenvalue moduli no more than $\gamma_j$.
Apparently, each $\im E_j$ (and thus  each $\im E_j\setminus \im E_{j+1}$) is invariant under the action of $A$ and all $A^k$ for $k\in\N$.
We get Theorem \ref{Nayak}(iii). 
Finally, every nonzero vector $x\in \C^n$ can be uniquely written as 
$x=x_1+\cdots+x_s$ where $x_i\in \im Q_i$ for $1\le i\le s$. Note that $\im Q_i$ is the $\gamma_i$-eigenspace for the hyperbolic element $H$. 
Given   $x$, suppose $x\in\im E_j\setminus \im E_{j+1}$ for certain $1\le j\le s$. Then
$x_1=\cdots=x_{j-1}=0$ and $x_j\ne 0$, so that
\begin{equation}
    A^m x=\sum_{i=j}^{s} A^m  x_i
    =\sum_{i=j}^{s} E^mU^m (H^m x_i)=\sum_{i=j}^{s} \gamma_i^{m} E^mU^m  x_i.
\end{equation}
Theorem \ref{thm: Yamamoto extension} shows that $\lim_{m\to\infty} |E^mU^m |^{1/m}=I_n$.
Remark \ref{rectangle B_mA_m limit} implies that for each nonzero $y\in\C^n$,  
\begin{equation}\label{E^mU^m y}
    \lim_{m\to\infty} \|E^mU^m y\|^{1/m}=\lim_{m\to\infty} \| y\|^{1/m}=1.
\end{equation}
Therefore,
\begin{eqnarray}
    \lim _{m\to \infty}\|A^m x\|^{1/m}
    &=& \lim _{m\to \infty}\|\sum_{i=j}^{s} \gamma_i^{m} E^mU^m  x_i \|^{1/m}
\notag\\
&=& \gamma_j \lim _{m\to \infty}\|E^mU^m x_j+\sum_{i=j+1}^{s} ( \gamma_i/\gamma_j)^{m} E^mU^m  x_i \|^{1/m}.
\notag
\end{eqnarray}
By \eqref{E^mU^m y}, $\lim _{m\to \infty}\|E^mU^m x_j\|^{1/m}=1$ and $\lim _{m\to \infty}( \gamma_i/\gamma_j)^{m} E^mU^m  x_i=0$ for each $j+1\le i\le s.$
Therefore, it is not hard to get Theorem \ref{Nayak}(ii):
\begin{equation}
    \lim _{m\to \infty}\|A^m x\|^{1/m}=\gamma_j.
\end{equation}
\end{rem}

The positive semidefinite part $|A|'$ of the other polar decomposition of $A = |A|'U$, where $U\in \U(n)$, is
 \begin{equation}
    |A|' := (AA^*)^{1/2}=|A^*|.
\end{equation}
It is easy to see that the hyperbolic component of $A^*$ is $H^*$ where $H$ is defined by \eqref{hyperbolic} for all $A\in M_n(\C)$. 
Moreover, if $A\in\GL_n(\C)$ and $A$ has the CMJD $A=EHU$, 
then $A^*=E^* H^* U^*$ is the CMJD of $A^*$. We have the following result.
\begin{cor}
\label{thm: Yamamoto extension 2}  
Suppose $A\in M_n(\C)$. Let $H=MDM^{-1}$ be the hyperbolic element of $A$ given  in \eqref{hyperbolic}, in which 
$D=\diag(|\lambda_1 (A)|,\ldots,|\lambda_n (A)|).$
For any $B, C\in\GL_n(\C)$,  let 
$BM=QR$  in which  $Q$ is   unitary and $R$ is upper  triangular. Then
\begin{equation}\label{BA^mC limit 2}
\lim_{m\to\infty} {|BA^mC|'}^{1/m} = QDQ^*,
\end{equation}
 in which the limit is independent of $C$  and the choice of $M$. 
\end{cor}
\begin{proof} We have
    $$\lim_{m\to\infty} {|BA^mC|'}^{1/m} = \lim_{m\to\infty} |C^*(A^*)^mB^*|^{1/m}.
$$
The hyperbolic component of $A^*$ is $H^*=M^{-*}DM^*$, and $M^*B^*=R^*Q^*$ where
$R^*$ is lower triangular and $Q^*$ is unitary. 
Then apply Theorem \ref{thm: Yamamoto extension} to get the limit.
\end{proof}

\section{Semisimple Lie Group Extensions}

In \cite{TamH06},  
Yamamoto's theorem \eqref{yamamoto} was extended in the context of  real  semisimple Lie groups. 
Here we will extend Theorem \ref{thm: Yamamoto extension} 
 in the same context.
Theorem \ref{thm: Yamamoto extension} involves  CMJD, polar decomposition/SVD, and the QR decomposition
$C^*M^{-*}=Q^*L^*$ in $\GL_n(\C)$.  They correspond to 
CMJD, Cartan decomposition/$KA_+K$ decomposition, and the Iwasawa decomposition in the real semisimple Lie groups. 

\subsection{Decompositions on real semisimple Lie groups}

Let $G$ be a  noncompact  connected real semisimple Lie group with the corresponding Lie algebra $\g$, which must be real semisimple. 
The Cartan decompositions on $\g$ and $G$ are discussed in \cite[VI.2, VI.3]{Knapp}. 
Explicitly, given a Cartan involution $\theta$ of $\g$, let $\k$ (resp. $\p$) be the $+1$ (resp. $-1$) eigenspace of $\theta$. The decomposition
\begin{equation}
   \g = \k \oplus\p 
\end{equation}
is a {\em Cartan decomposition of $\g$}.  
Let $\Theta$ be the global Cartan involution on $G$ with th differential $\theta$. 
Denote 
\begin{equation}
    g^*=\Theta(g^{-1}),\qquad g\in G.
\end{equation}
The analytic subgroup $K$ of $G$ with Lie algebra $\k$ is exactly the subgroup of $G$ fixed by $\Theta$. 
Denote 
\be{equation}{
P:=\exp \p.
}
For $p=\exp X\in P$ where $X\in\p$ and $r\in\R$, 
\begin{equation}
    p^r:=\exp(rX)\in P,
\end{equation}
which is well defined since
the map   $K\times \p\to G$ defined by $(k,X)\mapsto k\exp X$ is a diffeomorphism onto. 
Every $g\in G$ can be uniquely written as
\begin{equation}
    g=k(g) p(g), \qquad k(g)\in K,\ p(g)\in P,
\end{equation} 
 which is called the {\em Cartan decomposition of $g$ in $G$}. 
 We have $k^*=k^{-1}$ for $k\in K$ and $p^*=p$ for $p\in P$. Hence for $g\in G$, we have $g^*g\in P$ and
 \begin{equation}
     p(g)=(g^*g)^{1/2}.
 \end{equation}
For simplicity, we denote $|g|:=p(g)$, which is consistent with the notation in matrices.

The analogy of matrix QR decomposition to   semisimple Lie groups is the Iwasawa decomposition \cite[VI.4]{Knapp}.  
Let $\a$ be a maximal abelian subspace of $\p$. With respect to the $\ad$-action, $\g$ has the {\em restricted root space decomposition}
\begin{equation}
    \g=\g_0\oplus \bigoplus_{\lambda\in\Sigma} \g_{\lambda}
\end{equation}
 in which $g_{\lambda}$ and the set $\Sigma$ of restricted roots   are given by
\begin{eqnarray}
\g_{\lambda} &:=& \{X\in\g\mid (\ad H)X=\lambda(H)X \text{ for all } H\in \a\},
\\
\Sigma &:=& \{\lambda\in\a^*\setminus\{0\}\mid \g_{\lambda}\ne 0\}.
\end{eqnarray}
Fix a {\em closed} Weyl chamber $\a_+$ in $\a$. 
In Lie group $G$, set
\begin{equation}
A:=\exp \a,\qquad A_+ := \exp \a_+.
\end{equation}
The set $\Sigma^+$ of positive roots in the dual space $\a^*$ are also fixed by $\a_+$. Then
\begin{equation}
    \n := \bigoplus_{\lambda\in\Sigma^+} \g_{\lambda}
\end{equation}
is  a nilpotent Lie subalgebra of $\g$. {\em The Iwasawa decomposition of Lie algebra $\g$} is the vector space direct sum \cite[Proposition 6.43]{Knapp}:
\begin{equation}
    \g=\k\oplus\a\oplus\n.
\end{equation}
Let 
\begin{equation}
    N:=\exp\n.
\end{equation}
 The {\em Iwasawa decomposition}  of $G$:  
\begin{equation}
    K\times A\times N\to  G,\qquad (k, a, n)\mapsto kan. 
\end{equation}
is a diffeomorphism onto \cite[Proposition 6.46]{Knapp}.

 The  {\em $KA_+K$ decomposition}  of  $G$ says that \cite[Theorem 1.1]{Helgason}:
\begin{equation}
    G=KA_+K.
\end{equation}
In particular, we have $P=\Ad(K) A_+$ where $\Ad(k)a=kak^{-1}$. This decomposition is related to the Cartan decomposition in the way that
if $g=k_1 ak_2$ for $k_1,k_2\in K$ and $a\in A_+$, then $g=(k_1k_2) (k_2^{-1}ak_2)$ is the Cartan decomposition of $g$,
where $k_1k_2\in K$ and $k_2^{-1}ak_2\in P$.

    An element $h\in G$ is called {\em hyperbolic} if $h =
\exp X$ where $X\in \g$ is real semisimple, that is, $\ad X\in
\End(\g)$ is diagonalizable over $\R$. An element $u\in G$ is called
{\em unipotent} if $u = \exp(N)$ where $N\in \g$ is nilpotent, that
is, $\ad N\in \End(\g)$ is nilpotent.
 An element $e\in G$ is {\em elliptic} if $\Ad (e)\in \Aut(\g)$ is diagonalizable
over $\C$ with eigenvalues of modulus $1$. {\em The complete
multiplicative Jordan decomposition (CMJD)} 
for $G$ asserts that each $g\in G$ can be uniquely written as \cite [Proposition 2.1]{Kostant}
\begin{equation}
  g = e(g)h(g)u(g),
\end{equation}
where $e(g)$ is elliptic, $h(g)$ is hyperbolic, $u(g)$ is unipotent, and
the three elements  mutually commute. 
Moreover, an element $h\in G$ is
hyperbolic if and only if $h$ is conjugate to a unique element
$b(h)\in A_+$ \cite [Proposition 2.4]{Kostant}. Denote
\begin{equation}
    b(g):=b(h(g)).
\end{equation}

\subsection{The limit of $|g_1g^m g_2|^{1/m}$}

Theorem \ref{thm: Yamamoto extension} can be extended to real semisimple Lie groups  below.

\begin{thm}\label{thm: Yamamoto extension Lie group}
    Let $G$ be a  noncompact  connected real semisimple Lie group.   
   Then for any $g, g_1, g_2\in G$,  
\begin{equation}\label{lim power Lie group}
    \lim_{m\to\infty} |g_1 g^m g_2|^{1/m} = k b(g) k^*= k b(g) k^{-1},
\end{equation}
in which 
\begin{enumerate}
    \item $g=ehu$ is the CMJD of $g$ in $G$, where $e$ is elliptic, $h=qb(g)q^{-1}$ is hyperbolic  with $b(g)\in A_+$ and $q\in G$, and $u$ is unipotent, 
    \item $k$ comes from the Iwasawa decomposition $g_2^*q^{-*}=kan$, where $k\in K$, $a\in A$, and $n\in N $. 
\end{enumerate}
   In particular, the choices of $g_1$ and $q$ do  not affect the limit 
\eqref{lim power Lie group}.
\end{thm}

\begin{proof} Let $G$ have dimension $n$. 
We look at the adjoint representations  $\Ad: G\to \Aut(\g)$. 
There exists an orthonormal basis (with respect to the Killing form of $\g$) such that $\Ad (K)$ 
(resp. $\Ad (P)$, $\Ad (A)$, $\Ad (N)$) consists of unitary (resp. positive definite, positive diagonal, unit upper triangular) matrices in $\GL_n(\C)$. Moreover, we have $\Ad(g^*)=\Ad(g)^*$ and $\Ad(|g|) = |\Ad(g)|$ for $g\in G$ under the  basis. 

For $g\in G$ and the corresponding CMJD $g=ehu$, we have the decomposition
$\Ad (g)=\Ad (e) \Ad (h) \Ad (u)$ in which 
$\Ad(e)$ is diagonalizable
over $\C$ with eigenvalues of modulus $1$,
$\Ad(u)=\Ad(\exp(N))=\exp(\ad(N))$ is unipotent, and 
$\Ad(h)=\Ad(\exp(X))=\exp(\ad(X))$ is hyperbolic.
Thus 
 the CMJD of the matrix $\Ad (g)$ is $\Ad (g)=\Ad (e) \Ad (h) \Ad (u)$. 
 As $h$ is hyperbolic, there is $q\in G$ such that $h=qb(g)q^{-1}$,
  where $b(g)\in A_+$ \cite[Proposition 2.4]{Kostant}. So $\Ad(h)  = \Ad (q) \Ad (b(g)) \Ad(q)^{-1}$. Note that 
$(q^{-1}g_2 )^* =kan$ by Iwasawa decomposition. Thus 
$$
 \Ad (q^{-1}g_2 ) =\Ad (n^*)\Ad (a^*) \Ad(k^*) =
(\Ad (n))^* (\Ad (a)) (\Ad(k))^{-1}
$$ 
in which $(\Ad (n))^*$ is unit lower triangular and $\Ad (a)$ is diagonal. 
 By Theorem \ref{thm: Yamamoto extension}, 
$$
\lim_{m\to\infty}  |\Ad  (g_1)\Ad( g)^m \Ad(g_2) |^{1/m} = \Ad(k)^{-1} \Ad(b(g))\Ad(k) = \Ad(k^{-1}b(g)k),
$$ 
As the adjoint representation
$\Ad:  G\to \Aut \g$,  $g\mapsto \Ad(g)$ 
is continuous, we have
\begin{eqnarray*}
 \Ad(\lim_{m\to\infty}  |g_1 g^mg_2|^{1/m})&=&
\lim_{m\to\infty}  \Ad (|g_1 g^mg_2|^{1/m}) \\
&=&\lim_{m\to\infty}  (\Ad (|g_1 g^mg_2|))^{1/m} \\ &= &\lim_{m\to\infty}  |\Ad  (g_1 g^mg_2) |^{1/m} \\
&= & \lim_{m\to\infty}  |\Ad  (g_1)\Ad( g)^m \Ad(g_2) |^{1/m} \\
&=& \Ad(k^{-1}b(g)k).
\end{eqnarray*}
Hence 
every limit point of $\{|g_1 g^mg_2|^{1/m}\}_{m\in\N}$ has the form $zk^{-1}b(g)k$, where $z$ is in the center  $Z$ of $G$. 
On one hand, the CMJD of $zk^{-1}b(g)k$ is $z(k^{-1}b(g)k)$ in which $z$ is elliptic and $k^{-1}b(g)k$ is hyperbolic.  
On the other hand, every limit point of $\{|g_1 g^mg_2|^{1/m}\}_{m\in\N}\subseteq P$ must be in $P$, which is hyperbolic  \cite [Proposition 6.2]{Kostant}.  By the uniqueness of CMJD, we conclude that  $z$  is the identity. As a result, we have
$$
\lim_{m\to\infty}  |g_1 g^mg_2|^{1/m} = k^{-1}b(g)k.
$$
We are going to show that the limit \eqref{lim power Lie group} is independent of the choice of $q$.  Let $\hat q\in G$ such that $h=\hat qb(g)\hat q^{-1}$. Then $\hat q = qr$, where $r$ fixes $b(g)$ via conjugation, hat is, $rb(g)r^{-1} = b(g)$. Let
$g_2^*q^{-*}=kan$ and $g_2^*\hat q^{-*}=\hat k\hat a\hat n$ according to the Iwasawa decomposition.
By Theorem \ref{thm: Yamamoto extension}, $$\Ad(\lim_{m\to\infty}  |g_1 g^mg_2|^{1/m}) =  \Ad(k^{-1}b(g)k) = \Ad(\hat k^{-1}b(g)\hat k)$$ as the limit is independent of $\Ad(q)$. So $k^{-1}b(g)k =z (\hat k^{-1}b(g)\hat k)$, where $z\in Z$ is elliptic. As $k^{-1}b(g)k$ and $\hat k^{-1}b(g)\hat k$ are both in $P$ so they are hyperbolic. By the uniqueness of CMJD, $z$ is the identity and hence $k^{-1}b(g)k = \hat k^{-1}b(g)\hat k$. Thus the limit \eqref{lim power Lie group} is independent of the choice of $q$.

Similarly the limit is independent of $g_1$.
\end{proof}

\begin{rem}\rm 
Once we fix the Cartan decomposition, $G=KP$, the left side $ \lim_{m\to\infty} |g_1 g^m g_2|^{1/m}$ of \eqref{lim power Lie group} is clearly fixed. In other words, the right side $kb(g)k^{-1}$ of \eqref{lim power Lie group}  is independent of the choice of $A$ and $A_+$ which determines the fundamental roots and vice versa, and $N$, as long as the Iwasawa decomposition is of the form $KAN$. Here is an explanation. If we choose another maximal abelian subspace $\tilde \a$ of $\p$, set $\tilde A := \exp \tilde \a$  and  fix a positive Weyl chamber $\tilde A_+$, then \cite [p.378]{Knapp} $\tilde \a = \Ad (v) \a$   for some $v\in K$ and thus 
${\tilde A}=vAv^{-1}$; furthermore, we may choose $v$ in the way that ${\tilde A}_+ = vA_+v^ {-1}$ and thus $\tilde N = vNv^{-1}$.
So
 $$g_2^*\tilde q^{-*}= g_2^*(qv^{-1})^{-*} =g_2^*(q)^{-*}v^{-1}= kan v^{-1} = \tilde k \tilde a\tilde n,$$
 where
 $$
 \tilde k = kv^{-1}\in K,\quad  \tilde a = vav^{-1}\in \tilde A,\quad \tilde n = vnv^{-1}\in \tilde N.
 $$
Regarding the unique hyperbolic element $h$ of $g$,
$$
h  = qb(g)q^{-1} = { \tilde q} {\tilde b}(g) {\tilde q}^{-1},
$$
 where $ {\tilde b}(g)  = vb(g)v^{-1}\in \tilde A_+$ and ${ \tilde q} = qv^{-1}$. As a result, we have
 $$
 \tilde k\tilde b(g)\tilde k^{-1} =(kv^{-1})vb(g)v^{-1}(kv^{-1})^{-1}  =   kb(g)k^{-1},
 $$ 
 that is, the limit is independent the choice of $A$ and $A_+$, and thus independent of the choice of $A$ and $N$.
\end{rem}

\bibliography{Yamamoto}

\begin{thebibliography}{10}

\bibitem{S}
D.~Serre, {\em Matrices: theory and applications}, vol.~216 of {\em Graduate
  Texts in Mathematics}.
\newblock Springer, New York, second~ed., 2010.

\bibitem{HJ13}
R.~A. Horn and C.~R. Johnson, {\em Matrix analysis}.
\newblock Cambridge University Press, Cambridge, second~ed., 2013.

\bibitem{Y}
T.~Yamamoto, ``On the extreme values of the roots of matrices,'' {\em J. Math.
  Soc. Japan}, vol.~19, pp.~173--178, 1967.

\bibitem{TamH06}
T.-Y. Tam and H.~Huang, ``An extension of {Y}amamoto's theorem on the
  eigenvalues and singular values of a matrix,'' {\em J. Math. Soc. Japan},
  vol.~58, no.~4, pp.~1197--1202, 2006.

\bibitem{HuangTam08}
H.~Huang and T.-Y. Tam, ``Some asymptotic behaviors associated with matrix
  decompositions,'' {\em Int. J. Inf. Syst. Sci.}, vol.~4, no.~1, pp.~148--159,
  2008.

\bibitem{Nayak}
S.~Nayak, ``A stronger form of {Y}amamoto's theorem on singular values,'' {\em
  Linear Algebra Appl.}, vol.~679, pp.~231--245, 2023.

\bibitem{Bhatia07}
R.~Bhatia, {\em Positive definite matrices}.
\newblock Princeton Series in Applied Mathematics, Princeton University Press,
  Princeton, NJ, 2007.

\bibitem{Helgason}
S.~Helgason, {\em Differential geometry, {L}ie groups, and symmetric spaces},
  vol.~34 of {\em Graduate Studies in Mathematics}.
\newblock American Mathematical Society, Providence, RI, 2001.
\newblock Corrected reprint of the 1978 original.

\bibitem{Knapp}
A.~W. Knapp, {\em Lie groups beyond an introduction}, vol.~140 of {\em Progress
  in Mathematics}.
\newblock Birkh\"{a}user Boston, Inc., Boston, MA, second~ed., 2002.

\bibitem{Kostant}
B.~Kostant, ``On convexity, the {W}eyl group and the {I}wasawa decomposition,''
  {\em Ann. Sci. \'{E}cole Norm. Sup. (4)}, vol.~6, pp.~413--455 (1974), 1973.

\end{thebibliography}
\bibliographystyle{ieeetr}

\end{document}